\documentclass{article}
\usepackage{ijcai22}

\usepackage{multirow}   %导言区

\usepackage{times}
\usepackage{soul}
\usepackage{url}
\usepackage[hidelinks]{hyperref}
\usepackage[utf8]{inputenc}
\usepackage[small]{caption}
\usepackage{graphicx}
\usepackage{amsmath}
\usepackage{amssymb}
\usepackage{amsthm}
\usepackage{bm}
\usepackage{xcolor}
\usepackage{booktabs}
\usepackage{algorithm}
\usepackage{algorithmic}
\usepackage{lineno}
\newtheorem{example}{Example}
\newtheorem{theorem}{Theorem}
\newtheorem{proposition}{Proposition}
\newtheorem{assumption}{Assumption}
\newtheorem{problem}{Problem}
\newtheorem{definition}{Definition}
\newtheorem{remark}{Remark}

\title{Model Predictive Control with Reach-avoid Analysis}

\author{
Dejin Ren$^{1,2}$
\and
Wanli Lu$^3$\and
Jidong Lv $^3$\and
Lijun Zhang$^{1,2}$\And
Bai Xue$^{1,2}$\thanks{Corresponding author}
\affiliations
$^1$State Key Lab. of Computer Science, Institute of Software, CAS, Beijing, China\\
$^2$University of Chinese Academy of Sciences, Beijing, China\\
$^3$National Engineering Research Center of Rail Transportation Operation and Control System, Beijing Jiaotong University, Beijing, China \\
\emails
\{rendj, zhanglj, xuebai\}@ios.ac.cn,
\{luwanli, jdlv\}@bjtu.edu.cn
}

\begin{document}

\maketitle

\begin{abstract}
In this paper we investigate the optimal controller synthesis problem, so that the system under the controller can reach a specified target set while satisfying given constraints. Existing model predictive control (MPC) methods learn from a set of discrete states visited by previous (sub-)optimized trajectories and thus result in computationally expensive mixed-integer nonlinear optimization. In this paper a novel MPC method is proposed based on reach-avoid analysis to solve the controller synthesis problem iteratively. The reach-avoid analysis is concerned with computing a reach-avoid set which is a set of initial states such that the system can reach the target set successfully. It not only provides terminal constraints, which ensure feasibility of MPC, but also expands discrete states in existing methods into a continuous set (i.e., reach-avoid sets) and thus leads to nonlinear optimization which is more computationally tractable online due to the absence of integer variables. Finally, we evaluate the proposed method and make comparisons with state-of-the-art ones based on several examples. 
\end{abstract}

\section{Introduction}

Control synthesis is a fundamental problem which automatically constructs a control strategy that induces a system to exhibit a desired behavior. Due to the ability of handling control and state constraints and yielding high performance control systems \cite{camacho2013model}, control design methods based on MPC have gained great popularity and found wide acceptance in industrial applications, ranging from autonomous driving \cite{verschueren2014towards,kabzan2019learning} to large scale interconnected power systems \cite{ernst2008reinforcement,mohamed2011decentralized}.

In MPC design methods one issue is to guarantee feasibility of the successive optimization problems \cite{scokaert1999feasibility}. Because MPC is `greedy' in its nature, i.e., it only searches for the optimal strategy within a finite horizon, an MPC controller may steer the state to a region starting from where the violation of hard state constraints cannot be avoided. Although this feasibility issue could be solved by using a sufficiently long prediction horizon, we may not be able to afford the computation overhead due to the limited computational resources. Consequently, several solutions towards the feasibility issue are proposed \cite{zheng1995stability,zeng2021safety,ma2021feasibility}. 
Besides the feasibility issue of satisfying all hard state constraints \cite{zheng1995stability,zeng2021safety,ma2021feasibility}, a system is also desired to achieve certain performance objective in an optimal sense. In existing literature, stability performance of approaching an equilibrium within the MPC framework is extensively studied. One common solution of achieving stability is adding Lyapunov functions as terminal costs and/or corresponding invariant sets as terminal sets \cite{michalska1993robust,limon2005enlarging,limon2006stability,mhaskar2006stabilization,de2008lyapunov,wu2019control,grandia2020nonlinear}. This has motivated significant research work on applications of this control design to nonlinear processes \cite{yao2021state}. However, in many real applications stability performance is demanding. For example, a spacecraft rendezvous may require the chaser vehicle to be at a certain position relative to the target, moving towards it with a certain velocity. All of these quantities are specified with some tolerance, forming a target region in the state space. When the chaser enters that region a physical connection would be made and the maneuver is complete. However, since the target velocity is non-zero, the region is not invariant and stability cannot be achieved. Regarding this practical issue, the formulation in this paper replaces the notion of stability with reachability: given a target set, the system will achieve the reachability objective of reaching the target set in finite time successfully. To the best of our knowledge, the learning model predictive control (LMPC) proposed in \cite{rosolia2017learning}, which utilizes the historical data to improve suboptimal controllers iteratively, is the only method within the MPC framework which can solve this problem. However, it leads to mixed-integer nonlinear programming problems  which are fundamentally challenging to solve.

In this work, we consider control tasks where the goal is to steer the system from a starting configuration to a target set in finite time, while satisfying state and input constraints. The control task is formalized as a reach-avoid optimization problem. For solving the optimization problem, we propose a novel learning based MPC algorithm which fuses iterative learning control (ILC) \cite{arimoto1984bettering} and reach-avoid analysis in \cite{zhao2022inner}. The proposed algorithm utilizes MPC to iteratively improve upon a suboptimal controller. Based on a suboptimal controller obtained in the previous iteration, the reach-avoid analysis is to compute a set of initial states such that the closed-loop system can reach the target set without a violation of the state and input constraints. In our algorithm, the reach-avoid analysis not only provides terminal constraints which ensure feasibility of the MPC, but also expands the viability space of the system and thus facilitates exploring better trajectories with lower costs. Finally, several examples demonstrate the benefits of our algorithm over state-of-the-art ones.

The closest work to the present one is \cite{rosolia2017learning}. Due to the use of a set of discrete states visited by previous (sub-)optimized trajectories, computationally expensive mixed-integer nonlinear programming problems have to be addressed online in the proposed LMPC, limiting its application in practice. Our method is inspired by \cite{rosolia2017learning}. However, it incorporates the recently developed reach-avoid analysis technique and expands discrete states into a continuous (reach-avoid) set, which not only facilitates exploring better trajectories but also leads to a novel MPC formulation which involves solving more computationally tractable nonlinear programming problems online. Besides, this work is also close to the ones on model-based safe reinforcement learning, such as \cite{thomas2021safe,luo2021learning,wen2018constrained}. However, these studies prioritize enforcing the safety objective rather than achieving a joint objective of safety and reachability.

This paper is structured as follows. In Section \ref{sec:pre} the reach-avoid problem of interest is introduced. Then, we detail our algorithm for solving the reach-avoid problem in Section \ref{CS}, and evaluate it on several examples in Section \ref{sec:ex}. Finally, we conclude this paper in Section \ref{sec:con}.

\section{Preliminaries}
\label{sec:pre}
In this section we introduce the reach-avoid optimization problem of interest. Throughout this paper, $\mathbb{R}^n$ and $\mathbb{N}$ denote
the set of $n$-dimensional real vectors and non-negative integers, respectively. $\mathbb{R}[\cdot]$ denotes the ring of polynomials in variables given by the argument. $\sum [x]$ represents the set of sum-of-squares polynomials over variable $x$, i.e.,
$\sum [x] = \{p \in \mathbb{R}[x] \mid p = \sum_{i=1}^{k} q_i^{2}, q_{i} \in \mathbb{R}[x], i=1,2,\dots,k\}$.

\subsection{Problem Statement}
\label{sub:ps}
Consider a discrete-time dynamical system
\begin{equation}
\label{system}
\small
    x(t+1)=f(x(t),u(t)), t\in \mathbb{N},
\end{equation}
where $f(\cdot):\mathbb{R}^n \times \mathbb{R}^m \rightarrow \mathbb{R}^n$ is the system dynamics, $x(\cdot): \mathbb{N}\rightarrow \mathbb{R}^n$ is the state and $u(\cdot): \mathbb{N} \rightarrow \mathcal{U}\subseteq \mathbb{R}^m$ is the control.

\begin{definition}
    A control policy is a sequence of control inputs $\pi = \{u(i)\}_{i\in \mathbb{N}}$, where $u(i): \mathbb{N} \rightarrow \mathcal{U}$. Furthermore, we define $\Pi$ as the set of all control policies.
\end{definition}

Given a safe set $\mathcal{X}=\{x\in \mathbb{R}^n \mid w(x)\leq0\}$  and  a target set $\mathcal{T}=\{x\in \mathbb{R}^n\mid g(x)\leq0\}$ with $\mathcal{T}\subseteq \mathcal{X}$, a control policy $\pi\in \Pi$ is safe with respect to a state $x_0$, if the control policy $\pi$ will drive system \eqref{system} starting from $x(0)=x_0$ to reach the target set $\mathcal{T}$ in finite time while staying inside the safe set $\mathcal{X}$ before the first target hitting time. The associated cost of reaching the target set $\mathcal{T}$ safely is defined below.
\begin{definition}
Given a state $x(0)=x_0 \in \mathcal{X}\setminus \mathcal{T}$, and a safe policy $\pi=\{u(i)\}_{i\in \mathbb{N}}\in \Pi$, the cost with respect to the state $x_0$ and the policy $\pi$ is defined below:  
\begin{equation}
    \label{Q_fun}
    \small
    Q(x_0,\pi) =\sum_{k=0}^{L-1}h(x(k),u(k)), 
\end{equation}
 where $x(k+1)=f(x(k),u(k))$, $L=\min\{i \in \mathbb{N}\mid x(i)\in \mathcal{T}\}$, and 
 $h(\cdot,\cdot)$ is continuous and satisfies
    \begin{equation}
        \label{cost_fun_condition}
        \small 
        \begin{aligned}
         &h(x,u) \geq 0, \forall (x,u) \in \mathcal{X} \times  \mathcal{U}.
        \end{aligned}
    \end{equation}

\end{definition}

In this paper, we would like to synthesize 
a safe control policy $\pi \in \Pi$ with respect to a specified initial state $x_0$ such that system \eqref{system} will enter the target set $\mathcal{T}$ with the minimal cost in finite time while staying inside the safe set $\mathcal{X}$ before the first target hitting time. 
\begin{problem}
\label{prob}
We attempt to solve the following reach-avoid optimization problem:
\begin{equation}
\label{optimization1}
\small
    \begin{split}
    &J(x_0)=\min_{T,u(i),i=0,\ldots,T-1}\sum_{i=0}^{T-1} h(x(i),u(i))\\
    &\text{s.t.~}
    \begin{cases}
    &x(i+1)=f(x(i),u(i)), x(0)=x_0,\\
    &u(i)\in \mathcal{U}, x(i)\in \mathcal{X},\\
    &x(T)\in \mathcal{T},\\
    &i=0,\ldots,T-1,\\
    &T\in \mathbb{N},
    \end{cases}
\end{split}
\end{equation}
where $T$ is the first time of hitting the target set $\mathcal{T}$.
\end{problem}

Due to uncertainties on the target hitting time, it is challenging to solve optimization \eqref{optimization1} directly. However, the computational complexity can be reduced when searching for a feasible sub-optimal solution. In this paper we will adopt a learning strategy to solving \eqref{optimization1}, which iteratively improves upon already known suboptimal policies as the LMPC algorithm in \cite{rosolia2017learning}. Consequently, an initial feasible policy is needed, as formulated in Assumption \ref{ass_feasible_sys}.
\begin{assumption}
\label{ass_feasible_sys}
 Assume an initial policy $\pi^0=\{u^0(i)\}_{i\in \mathbb{N}}$, which can drive system \eqref{system} starting from $x_0$ to the target set $\mathcal{T}$in finite time safely, is available. The corresponding trajectory of system \eqref{system} can be obtained and is denoted by $\{x^0(i)\}_{0\leq i\leq L^0}$, where, 
 \[\small \begin{cases}
 &x^0(i+1)=f(x^0(i),u^0(i)), i=0,\ldots,L^0-1,\\ 
 &x^0(0)=x_0, x^0(L^0)\in \mathcal{T}.
 \end{cases}
 \]
\end{assumption}

\begin{remark}
The availability of a feasible control policy is not restrictive in practice for a number of applications. For instance, with race cars one can always run a path at very low speed to obtain a control policy.
\end{remark}

\subsection{Guidance-barrier Functions}
In this subsection we introduce guidance-barrier functions. They not only provide terminal constraints in our MPC method escorting system \eqref{system} to the target set $\mathcal{T}$ safely,  but also generate a set which curves out a viability space for system \eqref{system} to reach the target set $\mathcal{T}$ safely.

\begin{definition}
Given the safe set $\mathcal{X}$, target set $\mathcal{T}$ and a factor $\lambda \in (1,\infty)$, a bounded function $v(x): \mathcal{Y}\rightarrow \mathbb{R}$ is a guidance-barrier function of system \eqref{system} with the feedback controller $\hat{u}(\cdot):\mathcal{X}\rightarrow \mathcal{U}$, if it satisfies the following constraints:
\begin{equation}
\label{gbf_c}
\small
    \begin{cases}
    &v(f(x,\hat{u}(x)))\geq \lambda v(x), \forall x\in \mathcal{X}\setminus \mathcal{T},\\
    & v(x)\leq 0, \forall x\in \mathcal{Y}\setminus \mathcal{X},\\
    &v(x)\leq M, \forall x\in \mathcal{T},\\
         &v(x_0)>0,
    \end{cases}
\end{equation}
where $\mathcal{Y}=\{y\mid y=f(x,u),u\in \mathcal{U},x\in \mathcal{X}\}\cup \mathcal{X}$, and $M$ is a user-defined positive number.
\end{definition}
 When $f(x,\hat{u}(x))$ is polynomial over $x$, and $\mathcal{X}$ is semi-algebraic set, i.e., $f(x,\hat{u}(x)), w(x) \in \mathbb{R}[x]$, a set $\mathcal{Y}$ of the form $\{x\mid w_0(x) \leq 0\}$ with $w_0(x) \in \mathbb{R}[x]$ can be obtained using program (3) in \cite{zhao2022inner}.
\begin{remark}
\label{finite_horizon}
It is worth remarking here that if  \eqref{gbf_c} holds, for system \eqref{system} with $\hat{u}(\cdot):\mathcal{X}\rightarrow \mathcal{U}$, the induced trajectory starting from $x_0$ will hit the target set $\mathcal{T}$ within a bounded amount of time being less than or equal to $\log_{\lambda} \frac{M}{v(x_0)}$ (It can obtained according to $\lambda^T v(x_0)\leq v(x(T)) \leq M$, where $T$ is the first hitting time of $\mathcal{T}$.).
\end{remark}

A reach-avoid set $\mathcal{R}=\{x\in \mathcal{X}\mid v(x)>0\}$ can be computed via solving constraint \eqref{gbf_c}, which is a set of states such that there exists a control policy $\pi\in \Pi$ driving system \eqref{system} to enter the target set $\mathcal{T}$ in finite time while staying inside the safe set $\mathcal{X}$ before the first target hitting time.

\begin{theorem}
\label{reach-avoid-set}
Given the safe set $\mathcal{X}$, target set $\mathcal{T}$ and a factor $\lambda \in (1,\infty)$, if $v(x): \mathcal{Y}\rightarrow \mathbb{R}$ is a guidance-barrier function of system \eqref{system} with the feedback  controller $\hat{u}(\cdot):\mathcal{X}\rightarrow \mathcal{U}$, then $\mathcal{R}=\{x\in \mathcal{X}\mid v(x)>0\}$ is a reach-avoid set. 
\end{theorem}

Theorem \ref{reach-avoid-set} can be assured by Corollary 1 in \cite{zhao2022inner}. Due to space limitations, we omitted the proof herein. 

\begin{remark}
One of admissible  control policies $\pi\in \Pi$ such that system \eqref{system} satisfies the reach-avoid specification can be constructed by the feedback controller $\hat{u}(\cdot): \mathcal{X}\rightarrow \mathcal{U}$: when system \eqref{system} is in state $x(i) \in \mathcal{R}$, the corresponding control action is $u(i)=\hat{u}(x(i))$. We denote such a control policy by $\pi_{\hat{u}}$ in the rest of this paper. 
\end{remark}

\section{Reach-avoid Model Predictive Control}
\label{CS}
In this section we elucidate our learning-based algorithm for solving optimization \eqref{optimization1} in Problem \ref{prob}, which is built upon a so-called reach-avoid model predictive control (RAMPC). The proposed RAMPC is constructed based on a guidance-barrier function.

Our RAMPC algorithm is iterative and at each iteration it mainly consists of three steps. The first step is to synthesize a feedback controller by interpolating the suboptimal state-control pair obtained in the previous iteration. Then, a guidance-barrier function satisfying \eqref{gbf_c} with the synthesized feedback controller is computed. Finally, based on the computed guidance-barrier function a MPC controller, together with its resulting state trajectory, is generated online. The framework of the algorithm is summarized in Alg. \ref{RAMPC_framework}.
\begin{algorithm}
\caption{The framework for solving optimization \eqref{optimization1}.}
%\small
    \begin{algorithmic}
        \REQUIRE system \eqref{system}; initial state $x_0$; safe set $\mathcal{X}$; target set $\mathcal{T}$; control input set $\mathcal{U}$; feasible state-control trajectory
        $\{(x^0(i),u^0(i))\}_{i=0}^{L^0-1}$, of which the cost is $J^0(x_0)=\sum_{i=0}^{L^0-1}h(x^0(i),u^0(i))$; factor $\lambda$ and bound $M$ in \eqref{gbf_c}; maximum iteration number $K$; prediction horizon $N$ in RAMPC; termination threshold $\xi>0$. 
        \ENSURE Return $J^*(x_0)$.
        \FOR{$j=0:K$}
        \item[1] apply interpolation techniques to compute a feedback controller $\hat{u}^{j}(\cdot):\mathcal{X}\rightarrow \mathcal{U}$ based on the $j$-th state-control trajectory $\{(x^{j}(i),u^{j}(i))\}_{i=0}^{L^{j}-1}$;
        \item[2] compute $v^{j}(x)$ via solving \eqref{gbf_c} with $\hat{u}(x)=\hat{u}^j(x)$;
        \item[3] solve RAMPC optimization, which is constructed with $v^{j}(x)$, to obtain a state-control pair $\{(x^{j+1}(i),u^{j+1}(i))\}_{i=0}^{L^{j+1}-1}$ and the cost $J^{j+1}(x_0)$:
        \IF{$J^{j+1}(x_0)-J^{j}(x_0)\geq -\xi$}
        \STATE Return $J^*(x_0)=J^{j+1}(x_0)$;
        \ENDIF
        \ENDFOR
    \end{algorithmic}
      \label{RAMPC_framework}
\end{algorithm}

For solving \eqref{gbf_c} in Alg. \ref{RAMPC_framework}, when $f(x,\hat{u}(x))$ is polynomial over $x$, and $\mathcal{Y}$, $\mathcal{X}$ and $\mathcal{T}$ are semi-algebraic sets, i.e., $f(x,\hat{u}(x)), w_0(x), w(x), g(x) \in \mathbb{R}[x]$, the problem of solving constraint \eqref{gbf_c} can be transformed into a semi-definite programming problem \eqref{gbf_sdp},
\begin{equation}
\label{gbf_sdp}
\small
    \begin{cases}
    &v(f(x,\hat{u}(x))) - \lambda v(x)+ s_1(x)w(x) - s_2(x)g(x) \in \sum[x]\\
    & - v(x)+ s_3(x)w_0(x) - s_4(x)w(x) \in \sum[x]\\
    &M - v(x) + s_5(x)g(x)  \in \sum[x]\\
         &v(x_0)>0,%\\
%    &\textcolor{blue}{v(x(i))>0, i=0,\ldots,L-1,}
    \end{cases}
\end{equation}
where $s_j(x) \in \sum[x], j=1, \ldots, 5$, $v(x) \in \mathbb{R}[x]$.

Otherwise, sample-based approaches in the context of randomized algorithms \cite{tempo2013randomized} for robust convex optimization can be employed to solve constraint \eqref{gbf_c}. The basis of these approaches is the Almost Lyapunov condition \cite{liu2020almost}, which allow the Lyapunov conditions to be violated in restricted subsets of the space while still ensuring stability properties. Although results from these approaches lack rigorous guarantees, nice empirical performances demonstrated the practicality of these approaches in existing literature (e.g., \cite{chang2021stabilizing}). We will revisit it in our future work.

\begin{remark}
In Alg. \ref{RAMPC_framework}, the computations in the first two steps, i.e., synthesizing feedback controllers and solving \eqref{gbf_c}, can be carried out offline. Online computations occur only in the third step, which generate MPC controllers.
\end{remark}

In the following subsection we will introduce the RAMPC optimization in Alg. \ref{RAMPC_framework}.

\subsection{Reach-avoid Model Predictive Control}
\label{sub:rampc}
This subsection introduces the RAMPC optimization in Alg. \ref{RAMPC_framework}, which involves solving the MPC optimization problem online : 
\[J_{t \rightarrow t+N}^{\text{RAMPC},j}(x_t^j)=\min_{u_{k\mid t}}\Bigg[\sum_{k=0}^{N-1}h(x_{k\mid t}^j,u_{k\mid t}^j)+Q^{j}(x_{N\mid t}^j,\pi_{\hat{u}^{j-1}})\Bigg]\]
\begin{equation}
\label{MPC10}
\small
\begin{split}
    &\text{s.t.~}
    \begin{cases}
    &x_{{k+1\mid t}}^j=f(x_{k\mid t}^j,u_{k\mid t}^j),\\
    &u_{k\mid t}^j\in \mathcal{U}, x_{k\mid t}^j\in \mathcal{X},\\
    &k=0,1,\ldots,N-1, \\
    & v^{j-1}(x_{N\mid t}^j)\geq   \begin{cases}
    &\lambda^N v^{j-1}(x_0),\text{if~} t=0,\\
    &\lambda v^{j-1}(x_{N\mid t-1}^j), \text{otherwise.}
    \end{cases}\\ 
    &x_{0\mid t}^j= x^j_t 
    \end{cases}
\end{split}
\end{equation}
where the superscript $j$ is the iteration index, $x_{{k\mid t}}^j$ is the state predicted $k$ time steps ahead, computed at time $t$, initialized at $x_{0\mid t}^j=x_t^j$ with $x_0^j=x_0$, and similarly for $u_{k\mid t}^j$, and $Q^{j}(x_{N\mid t}^j)=\sum_{i=0}^{L-1}h(x_{i+N\mid t}^j,\hat{u}^{j-1}(x_{i+N\mid t}^j))$ is the cost with respect to  the state $x_{N\mid t}^j$ and the control policy $\pi_{\hat{u}^{j-1}}$.

\begin{figure}[htbp]
\center
\includegraphics[width=1.2in]{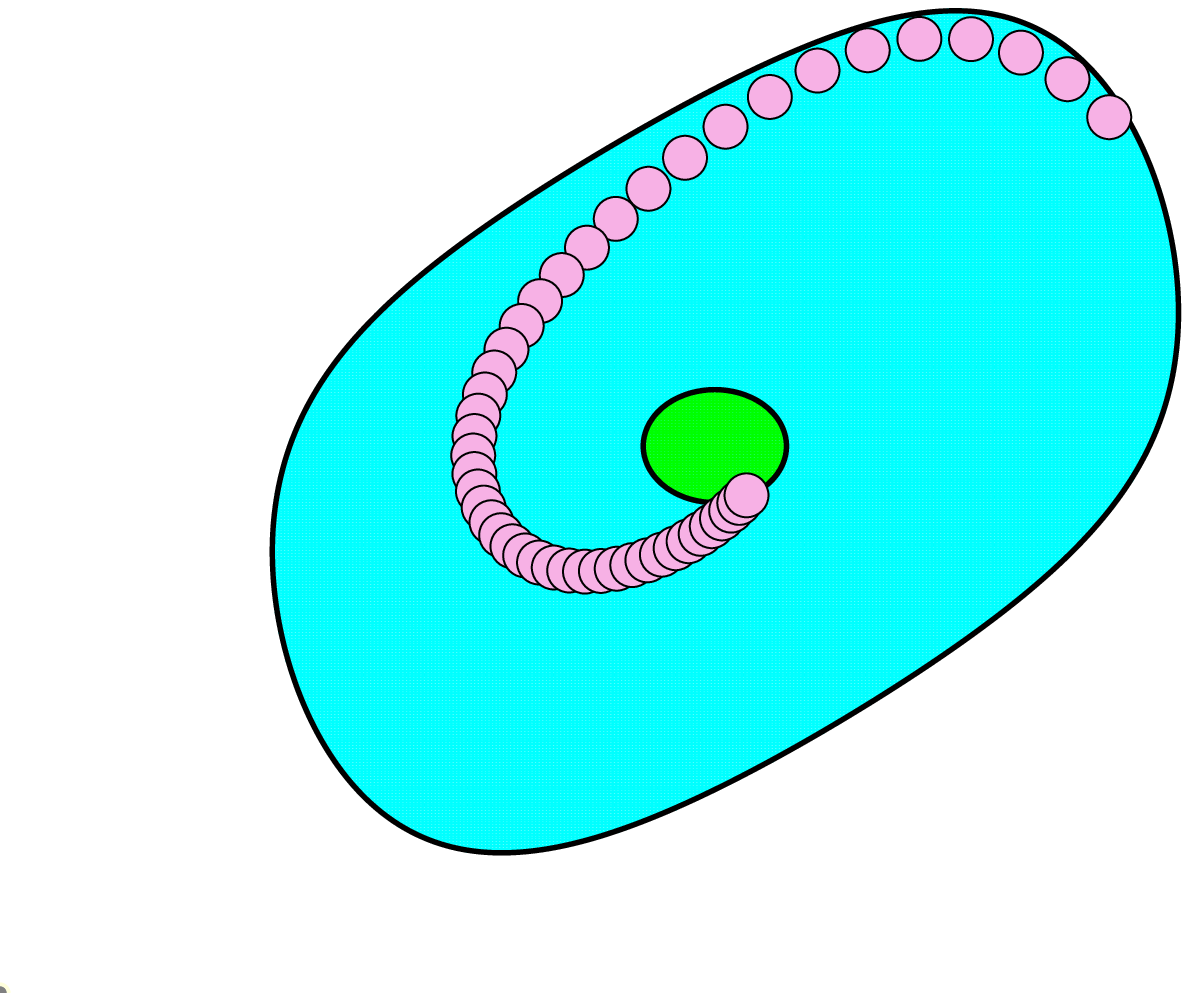}
\caption{An illustration of continuous (reach-avoid) sets in our RAMPC method and discrete states in the LMPC method, the green region denotes a target set.}
\label{ex_explain_fig}
\end{figure}

In \eqref{MPC10}, the terminal constraint 
\begin{equation}
\label{terminate_constraint}
\small 
    v^{j-1}(x_{N\mid t}^j)\geq    \begin{cases}
    &\lambda^N v^{j-1}(x_0),\text{if~} t=0,\\
    &\lambda v^{j-1}(x_{N\mid t-1}^j), \text{otherwise.}
    \end{cases}
\end{equation}
guarantees that the terminal state $x^j_{N\mid t}$ lies in the reach-avoid set $\mathcal{R}^{j-1}=\{x\in \mathcal{X}\mid v^{j-1}(x)>0\}$, which carves out a larger continuous viability space (visualized as the cyan region in Fig. \ref{ex_explain_fig}) for system \eqref{system} to be capable of achieving the reach-avoid objective. This is different from the LMCP method in \cite{rosolia2017learning}, which restricts terminal states within previously explored discrete states (visualized as pink points in Fig. \ref{ex_explain_fig}) and thus leads to computationally demanding mixed-integer nonlinear optimization. As to the practical computations of the cost $Q^{j}(\cdot,\pi_{\hat{u}^{j-1}}):\mathcal{R}^{j-1}\rightarrow \mathbb{R}$, we will give a detailed introduction in Subsection \ref{sec:PAC}.

Let
\begin{equation}
\label{optimal_solution}
  \small
\mathbf{u}_{0: N \mid t}^{*, j}=\{u_{i \mid t}^{*, j}\}_{i=0}^{N-1}  \text{~and~} \mathbf{x}_{0: N \mid t}^{*, j}=\{x_{i \mid t}^{*, j}\}_{i=0}^{N}
\end{equation}
be the optimal solution to \eqref{MPC10} at time $t$ and $J_{t \rightarrow t+N}^{\text {RAMPC},j}(x_t^j)$ the corresponding optimal cost. Then, at time $t$ of the iteration $j$, the first element $u^{j}(t)=u_{0 \mid t}^{*, j}$ of $\mathbf{u}_{0: N \mid t}^{*, j}$ is applied to the system \eqref{system} and thus the state of system \eqref{system} turns into  
\begin{equation}
\label{x_deploy}
\small
    x_{t+1}^j=x_{1\mid t}^{*,j} = f(x^j(t),u^j(t)),
\end{equation}
where $x^j(t)=x_{0\mid t}^{*,j}$, which will be the used to update the initial state in \eqref{MPC10} for subsequent computations. Once there exist $t\in \mathbb{N}$ and $l\leq N$ s.t. $x^j_{l\mid t} \in \mathcal{T}$, we terminate computations in this iteration: from the state $x_{0\mid t}^{*,j}$ on, we will apply the control actions $\{u^{*,j}_{i\mid t}\}_{i=0}^{l-1}$ successively to system \eqref{system}.

Let $\pi^{j}=\{u^{j}(t)\}_{0\leq t\leq L^j-1}$ be the control policy applied to system \eqref{system} by solving optimization \eqref{MPC10}. The resulting state-control pair $\{x^j(i),u^{j}(i)\}_{0\leq i\leq L^j-1 }$, where $x^{j}(L^j)\in \mathcal{T}$, is obtained. We can conclude that the performance of $\pi^j$ is no worse than that of $\pi^{j-1}$, i.e., $J^{j-1}(x_0)=\sum_{i=0}^{L^{j-1}-1} h(x^{j-1}(i),u^{j-1}(i))\geq \sum_{i=0}^{L^j-1} h(x^j(i),u^{j}(i))=J^j(x_0)$. This conclusion can be justified by Theorem \ref{J_non_increasing}.

Under Assumption \ref{ass_feasible}, we in Theorem \ref{feasibility} show that the RAMPC \eqref{MPC10} is feasible and system \eqref{system} with controllers computed by solving \eqref{MPC10} satisfies the reach-avoid property, and in Theorem \ref{J_non_increasing} show that the $j$-th iteration cost $J^j(x_0)$ is non-increasing as $j$ increases in RAMPC \eqref{MPC10}. Their proofs can be found in the appendix.

\begin{assumption}
\label{ass_feasible}
At each iteration $0\leq j\leq K$ the computed reach-avoid set $\mathcal{R}^{j}=\{x\in \mathcal{X}\mid v^{j}(x)>0\}$ via solving \eqref{gbf_c} is non-empty. In addition, we assume that it includes the state trajectory $\{x^{j}(i)\}_{i=0}^{L^{j}-1}$\footnote{This assumption can be realized by interpolating the feedback controller $\hat{u}^j(\cdot):\mathcal{X}\rightarrow \mathcal{U}$ s.t. $\hat{u}^j(x^j(i))=u^j(i)$, $i=0,\ldots,L^j-1$. This requirement mainly serves for our theoretical analysis since it can guarantee that our algorithm can iteratively improve the policy (Theorem \ref{J_non_increasing}). However, in practical this requirement is not indispensable regarding efficient computations.}.
\end{assumption}

\begin{theorem}
\label{feasibility}
   For each iteration $j\leq K$, the RAMPC \eqref{MPC10} is feasible for $1\leq t \leq L^j-1$; also, system \eqref{system} with controllers obtained by solving RAMPC \eqref{MPC10} can reach the target set $\mathcal{T}$ within the time interval $[0,\log_{\lambda}\frac{M}{v_0}]$ while satisfying all of state and input constraints.
   \end{theorem}

\begin{theorem}
    \label{J_non_increasing}
    Consider system \eqref{system} with controllers obtained by solving \eqref{MPC10}. Then, the iteration cost $J^j(x_0)$ is non-increasing with the iteration index $j$.
\end{theorem}

\subsection{Estimating Terminal Costs via Scenario Optimization}
\label{sec:PAC}
In practice, the exact terminal cost  $Q^j(\cdot,\pi_{\hat{u}^{j-1}}): \mathcal{R}^{j-1}\rightarrow \mathbb{R}$ is challenging, even impossible to obtain. In this subsection, we present a linear programming method based on the scenario optimization \cite{calafiore2006scenario} to compute an approximation of the terminal cost  $Q^j(\cdot,\pi_{\hat{u}^{j-1}}): \mathcal{R}^{j-1}\rightarrow \mathbb{R}$ in the probably approximately correct (PAC) sense.

Let $\{x_{i}\}_{i=1}^{N'}$ be the independent samples extracted from $\mathcal{R}^{j-1}=\{x\in \mathcal{X}\mid v^{j-1}(x)>0\}$ according to the uniform distribution $\mathbb{P}$, and $Q^j(x_{i},\pi_{\hat{u}^{j-1}})$ be the corresponding cost of the roll-out from $x_{i}$ with the controller $\pi_{\hat{u}^{j-1}}$. Such cost can be simply computed by summing up the cost along the closed-loop realized trajectory until the target set $\mathcal{T}$ is reached. Finally, using the set of sample states and correspding costs $\{(x_i,Q^j(x_{i},\pi_{\hat{u}^{j-1}}))\}_{i=1}^{N'}$, we approximate the terminal cost $Q^j(\cdot,\pi_{\hat{u}^{j-1}}): \mathcal{R}^{j-1}\rightarrow \mathbb{R}$ using the scenario optimization \cite{calafiore2006scenario}.

A linearly parameterized model template $Q^j_a(c_1,\ldots,c_l,x)$, $k\geq 1$ is utilized, which is a linear function in unknown parameters $\bm{c}=(c_1,\ldots,c_l)$ but can be nonlinear over $x$. Then we construct the following linear program over $\bm{c}$ based on the family of given datum $\{(x_i,Q^j(x_{i},\pi_{\hat{u}^{j-1}}))\}_{i=1}^{N'}$:
\begin{equation}
\label{co}
\small
    \begin{split}
        &\min_{\bm{c},\delta}\delta\\
        \text{s.t.~}&
        \begin{cases}
        &Q^j_a(\bm{c},x_i)-Q^j(x_{i},\pi_{\hat{u}^{j-1}})\leq \delta,\\
        &Q^j(x_{i},\pi_{\hat{u}^{j-1}})-Q^j_a(\bm{c},x_i)\leq \delta,\\
        &-U_{c}\leq c_l \leq U_{c},\\
        &i=1,\ldots,N',\\
        &0\leq \delta,
        \end{cases}
    \end{split}
\end{equation}
where $U_c\geq 0$ is a pre-specified upper bound for $c_i$, $i=1,\ldots,l$.

Denote the optimal solution to \eqref{co} by $(\bm{c}^*,\delta^*)$. The discrepancy between $Q^j_a(\bm{c}^*,x)$ and $Q^j(x,\pi_{\hat{u}^{j-1}})$ is formally characterized by two parameters: the error probability $\epsilon \in (0,1)$ and confidence level $\beta\in (0,1)$. 

\begin{theorem}
\label{sample_number}
    Let $(\bm{c}^*,\delta^*)$ be an optimal solution to \eqref{co}, $\epsilon\in (0,1)$, $\beta\in (0,1)$ and 
    \begin{equation*}
    \small
    \epsilon\geq \frac{2}{N'}(\ln \frac{1}{\beta}+l+1).
    \end{equation*}
    Then we have that with at least $1-\beta$ confidence, 
    \[\small \mathbb{P}(\{x\in \mathcal{R}\mid |Q^j_a(\bm{c}^*,x)-Q^j(x,\pi_{\hat{u}^{j-1}})|\leq \delta^*\})\geq 1-\epsilon.\]
\end{theorem}

Thus, we relax optimization \eqref{optimization1} into the following form:
\[\small J_{a,t \rightarrow t+N}^{\text {RAMPC}, j}(x_t^j)=\min_{u_{k\mid t}}\Bigg[\sum_{k=0}^{N-1}h(x_{k\mid t},u_{k\mid t})+Q^{j-1}_a(\bm{c}^*,x_{N\mid t})\Bigg]\]
\begin{equation}
\label{MPC10a}
\begin{split}
    \small
    &\text{s.t.~}
    \begin{cases}
    &x_{{k+1\mid t}}=f(x_{k\mid t},u_{k\mid t}),\\
    &u_{k\mid t}\in \mathcal{U}, x_{k\mid t}\in \mathcal{X}, \\
    &k=0,1,\ldots,N-1, \\
    & v(x_{N\mid t})\geq    \begin{cases}
    &\lambda^N v(x_0),\text{if~} t=0,\\
    &\lambda v(x_{N\mid t-1}), \text{otherwise.}
    \end{cases}\\ 
    &x_{0\mid t} = x^j_t 
    \end{cases}
\end{split}
\end{equation}
where $Q^{j-1}_a(\bm{c}^*,x_{N\mid t})$ is the approximate terminal cost at the state $x_{N\mid t}$, which is obtained via solving \eqref{co}.

\begin{proposition}
\label{probability theorem}
   With at least $1-\beta$ confidence,
   \begin{equation*}
   \small
   |J_{a,t \rightarrow t+N}^{\text {RAMPC}, j}(x_t^j)-J_{t \rightarrow t+N}^{\text {RAMPC}, j}(x_t^j)|\leq \delta^*
   \end{equation*}
holds with the probability larger than or equal to $1-\epsilon$, i.e., 
  \begin{equation*}
   \small
   \mathbb{P}(\{x_t^j\in \mathcal{R}^{j-1}\mid |J_{a,t \rightarrow t+N}^{\text {RAMPC}, j}(x_t^j)-J_{t \rightarrow t+N}^{\text {RAMPC}, j}(x_t^j)|\leq \delta^*\})\geq 1-\epsilon.
   \end{equation*}
\end{proposition}

Relying on \eqref{MPC10a}, we summarize our RAMPC algorithm for solving optimization \eqref{optimization1} in Algorithm \ref{RAMPC_alg}.
\begin{algorithm}
\caption{The RAMPC algorithm for solving \eqref{optimization1}.}
    \begin{algorithmic}
        \REQUIRE system \eqref{system} with an initial state $x_0$, a safe set $\mathcal{X}$, a target set $\mathcal{T}$ a control input set $\mathcal{U}$ and a feasible state-control trajectory
        $\{(x^0(i),u^0(i))\}_{i=0}^{L^0-1}$, of which the corresponding cost is $J^0(x_0)=\sum_{i=0}^{L^0-1}h(x^0(i),u^0(i))$; factor $\lambda$ and bound $M$ in \eqref{gbf_c}; iteration number $K$, prediction horizon $N$ and termination threshold $\xi>0$; probability error $\epsilon$ and confidence level $\delta$ in PAC approximations. 
        \ENSURE Return $J^*(x_0)$.
        \FORALL{$j=0:K$}
        \item[1] apply interpolation techniques to compute $\hat{u}^{j}(\cdot):\mathcal{X}\rightarrow \mathcal{U}$ for the $j$-th state-control trajectory $\{(x^{j}_i,u^{j}_i)\}_{i=0}^{L^{1}-1}$;
        \item[2] obtain $\mathcal{R}^{j}=\{x\in \mathcal{X}\mid v^{j}(x)>0\}$ via solving \eqref{gbf_c} with $\hat{u}(x)=\hat{u}^j(x)$;
        \item[3] compute the PAC terminal cost $Q^{j}_a(\cdot):\mathcal{R}^{j}\rightarrow \mathbb{R}$ via solving optimization \eqref{co};
         \item[4] Solving MPC optimization \eqref{MPC10a} with $v^j$ to obtain a state-control pair $\{(x^{j+1}(i),u^{j+1}(i))\}_{i=0}^{L^{j+1}-1}$ and the cost $J^{j+1}(x_0)$:
        \IF{$|J^{j+1}(x_0)-J^{j}(x_0)|\leq \xi$}
        \STATE Return $J^*(x_0)=J^{j+1}(x_0)$;
        \ENDIF   
         \ENDFOR
    \end{algorithmic}
    \label{RAMPC_alg}
\end{algorithm}

\section{Examples}
\label{sec:ex}
In this section we evaluate our RAMPC algorithm, i.e., Alg. \ref{RAMPC_alg}, and make comparisons with the LMPC algorithm in \cite{rosolia2017learning} on several examples. All the experiments were run on MATLAB 2022b with CPU 12th Gen Intel(R) Core(TM) i9-12900K
 and RAM 64 GB. Constraint \eqref{gbf_c} is solved by encoding it into sum-of-squares constraints which is treated by the semi-definite programming solver MOSEK; the nonlinear programming \eqref{MPC10} and the mixed-integer nonlinear programming in the LMPC algorithm in \cite{rosolia2017learning} are solved using YALMIP \cite{lofberg2004yalmip}. In addition, we take $\hat{u}^{j}(\cdot)$ as a linear function to interpolate the $j$-th state-control trajectory $\{(x^{j}_i,u^{j}_i)\}_{i=0}^{L^{j}-1}$ at each iteration $j \geq 1$. The configuration parameters in Alg. \ref{RAMPC_alg} for all examples are shown in Table \ref{para_table}.

\begin{table}[htbp]
\small
    \centering
    \begin{tabular}{|c|c|c|c|c|c|c|c|}
    \hline
         Example&$\lambda$& $M$&$N$&$K$&$\xi$&$\delta$&$\epsilon$ \\\hline
          Ex. \ref{ex1}&1.001&1&4&8&0.1&0.1&0.1 \\\hline
            Ex. \ref{ex3}&1.001&1&3&8&0.1 &0.05&0.05 \\\hline
             Ex. \ref{ex4}&1.001&1&4&6&0.002 &0.1&0.1 \\\hline
    \end{tabular}
    \caption{Configuration parameters in Alg. \ref{RAMPC_alg} for examples.}
    \label{para_table}
\end{table}

\begin{example}

\label{ex1}
Consider the drone system from \cite{rosolia2017learning},

\begin{equation*}
\label{ex1system}
\small 
x_{t+1}=\left[\begin{array}{ll}1 & dt \\ 0 & 1\end{array}\right] x_t+\left[\begin{array}{l}0 \\ 1\end{array}\right] u_t
\end{equation*}
where $x_t = (p_t, v_t)^\top$, $p_t$ and $v_t$ are respectively the position and velocity of the drone at time $t$, $u_t$ is the control input and $dt$ is the control interval which is equal to $0.1$.

We assume $\mathcal{X}=\{(p,v)^{\top}\mid  \frac{p^2}{8^2}+\frac{v^2}{8^2}-1\leq 0\}$, initial state $x_0 = (4,-6)^\top$, target set $\mathcal{T}=\{(p,v)^{\top}\mid p^2+v^2-0.5^2\leq 0\}$ and control input set $\mathcal{U} = \{u \mid -0.5 \leq u \leq 0.5\}$. The set $\mathcal{Y}=\{(p,v)^{\top}\mid  \frac{p^2}{8^2}+\frac{v^2}{8^2}-2\leq 0\}$ in constraint \eqref{gbf_c} is obtained by solving a semi-define program as in \cite{zhao2022inner}.  The cost in optimization \eqref{optimization1} is $h\left(x, u\right)=\left\|x\right\|_2^2+\left\|u\right\|_2^2$. For simplicity of presentation herein, we do not present the initial state-control trajectory which is a long sequence. The initial state-control trajectory is induced by the feedback controller $\hat{u}^0(p,v)=-0.04p-0.1v$.  For scenario optimization \eqref{co}, we use a template of the polynomial form $Q^j_a(\bm{c},p,v) = c_{1}+c_2p+c_3v+c_4pv+c_5p^2+c_6v^2$ and $N'=207$ sampled points which can be computed via Theorem \ref{sample_number}. 

We compare our RAMPC algorithm with the LMPC one in \cite{rosolia2017learning} with the same prediction horizon and termination conditions (i.e., $N$ and $\xi$). Since this system is linear, instead of using a set of explored discrete states, the terminal constraint in the LMPC algorithm can be constructed by its convex hull, thus resulting in non-linear programs instead of mixed-integer nonlinear programs as mentioned in \cite{rosolia2017learning}. The performances of trajectories (shown in Figure \ref{ex1_figure}) generated by both are compared. The iteration costs and computation times at each iteration of RAMPC and LMPC are presented in Table \ref{ex1_table}. It is observed that the iteration costs in our algorithm drop more quickly than those in the LMPC one. Also, our RAMPC algorithm is more efficient than the LMPC one. RAMPC terminates after $2$ iterations with the computation time of $6.3304s$, but LMPC converges more slowly and has to take $5$ iterations with the computation time of $8.2259s$.
\begin{table}[tbp]
    \centering
    \setlength{\tabcolsep}{1.5mm}{
    \begin{tabular}{*{5}{c}}
  \toprule
   \multirow{2}*{\textbf{Iteration}}&
   \multicolumn{2}{c}{\textbf{Iteration Cost}} & \multicolumn{2}{c}{\textbf{Time Cost}(seconds)}\\
  \cmidrule(lr){2-3}  \cmidrule(lr){4-5}
  &RAMPC & LMPC &RAMPC & LMPC \\
  \midrule
   0 & 369.8267 & 369.8267\\
   1 & 215.1007  & 222.2482& 3.3204 & 3.1741 \\
   2 & 215.1002 & 217.3604& 3.0100 & 1.5634 \\
   3 & - & 215.3008& - & 1.1980 \\
   4 &   -      & 215.1003  & - & 1.1023  \\
   5 & - & 215.1044 & - & 1.1881  \\
   \textbf{total time} & & & $\bm{6.3304}$ & 8.2259  \\
  \bottomrule
\end{tabular}}
    \caption{Iteration cost and computation time for Example \ref{ex1}.}
    \label{ex1_table}
\end{table}

\begin{figure}[htbp]
\center
\includegraphics[width=2.5in]{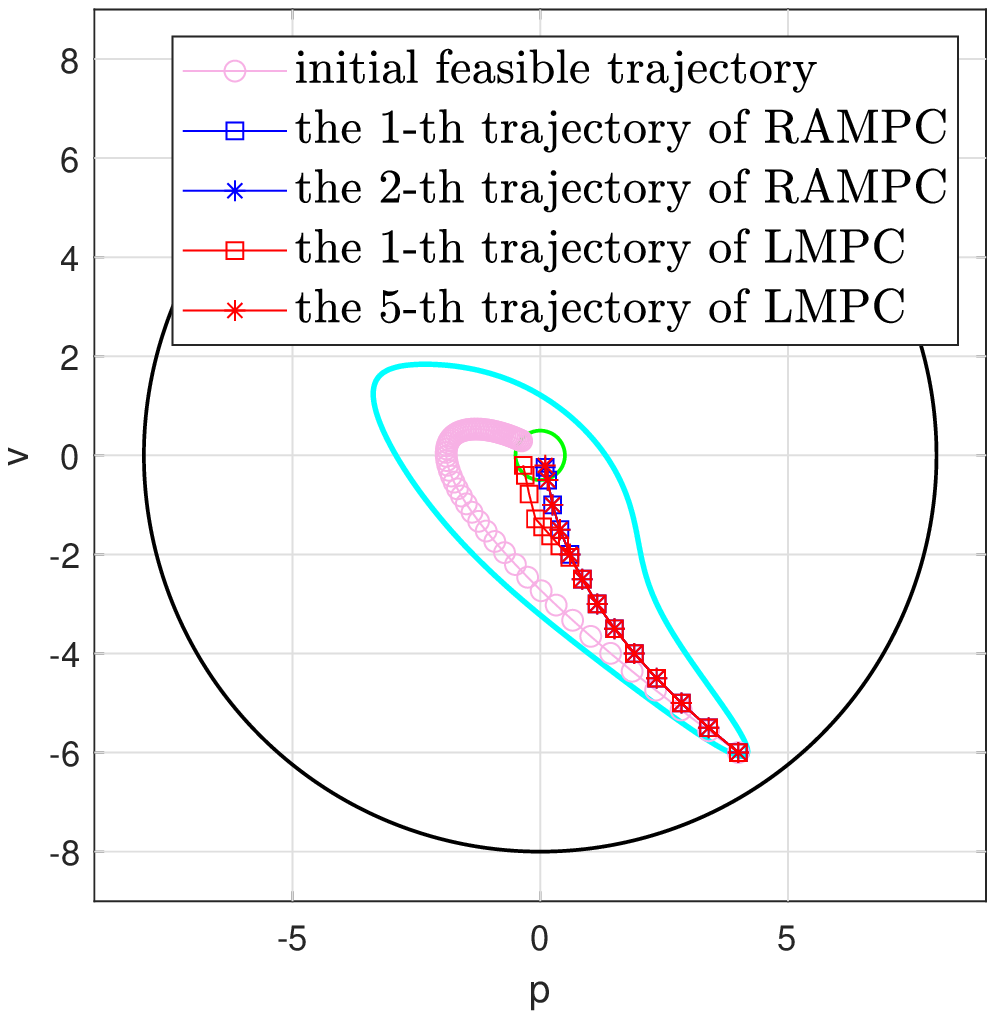} 
\caption{Trajectories in Example \ref{ex1} (\textcolor{green}{Green}, black and \textcolor{cyan}{cyan} curve- \textcolor{green}{$\partial \mathcal{T}$},  \textcolor{black}{$\partial \mathcal{X}$} and \textcolor{cyan}{$\partial \mathcal{R}^0$}).}
\label{ex1_figure}
\end{figure}

\end{example}

\begin{example}
\label{ex3}
Consider the Euler version with the time step $\Delta t = 0.05$ of the controlled reversed-time Van der Pol oscillator in \cite{drazin1992nonlinear}:

\begin{equation*}
\small
\begin{cases}
&x_1(t+1)=x_1(t)-\Delta t x_2(t) \\
&x_2(t+1)=x_2(t)-\Delta t ((1-x_1^2(t)) x_2(t)-x_1(t))+u(t),\\  
\end{cases}
\end{equation*}
with the safe set $\mathcal{X}=\{(x_1,x_2)^{\top}\mid  \left(\frac{x_1}{2}\right)^2+\left(\frac{x_2}{2}\right)^2-1\leq 0\}$, initial state $x_0 = (1.2,1)^\top$, target set $\mathcal{T}=\{(x_1,x_2)^{\top}\mid x_1^2+x_2^2-0.2^2\leq 0\}$ and input set $\mathcal{U} = \{u \mid -0.5 \leq u \leq 0.5\}$. 

In \eqref{gbf_c}, the set 
$\mathcal{Y}=\{(x_1,x_2)^{\top}\mid  \left(\frac{x_1}{2}\right)^2+\left(\frac{x_2}{2}\right)^2-2\leq 0\}$ is computed by solving a semi-define program in \cite{zhao2022inner}. In addition,  the cost $h\left(x,u\right)=\left\|x\right\|_2^2+\left\|u\right\|_2^2$, where $x = (x_1,x_2)^\top$, is adopted. The initial trajectory is induced by the  controller $\hat{u}^0(\bm{x})\equiv 0$. For scenario optimization \eqref{co}, we use a  template of the polynomial form $Q^j_a(\bm{c},x) = c_{1}+c_2x_1+c_3x_2+c_4x_1x_2+c_5x_1^2+c_6x_2^2$ and $N'=428$ samples which can be computed via Theorem \ref{sample_number}. LMPC uses the same prediction horizon and termination conditions (i.e., $N$ and $\xi$). 

Some trajectories generated  by our RAMPC algorithm and the LMPC algorithm are visualized in Figure \ref{ex3_figure}. Our RAMPC algorithm terminates after the third iteration, but the LMPC one does not terminate after eight iterations. Figure \ref{ex3_figure} shows that the initial trajectory and the trajectory generated by the first iteration in the LMPC algorithm are too close to be indistinguishable within the first few time steps, which on the other hand further reflects that the controller generated in the first iteration of the LMPC algorithm may not improve the performance induced by the initial control policy. Besides, it is observed that the eighth trajectory from the LMPC algorithm looks not stable and has strong oscillations, which are not expected in practice. In contrast, the trajectory generated by our algorithm is smoother.

Table \ref{ex3_it_cost} summarizes the iteration costs and the computation times at each iteration of our RAMPC algorithm and the LMPC one. The iteration cost induced by the initial trajectory is $64.3087$, after three iterations our RAMPC algorithm reduces the cost to $29.2824$ with the computation time of $31.2841s$. In contrast, the LMPC algorithm needs more than 8 iterations with the computation time of almost more than $3$ hours to achieve the same cost. This striking contrast definitely demonstrates that our RAMPC algorithm can solve optimization \eqref{optimization1} more efficiently for some cases.

\begin{table}[htbp] 
    \centering
    \setlength{\tabcolsep}{1.5mm}{
     \begin{tabular}{*{5}{c}}
  \toprule
   \multirow{2}*{\textbf{Iteration}}&
   \multicolumn{2}{c}{\textbf{Iteration Cost}} & \multicolumn{2}{c}{\textbf{Time Cost}(seconds)}\\
  \cmidrule(lr){2-3}  \cmidrule(lr){4-5}
  &RAMPC & LMPC &RAMPC & LMPC \\
  \midrule
   0 & 64.3087 & 64.3087 \\
   1 & 30.9693 & 40.0714& 10.6583 & 45.0 \\
   2 & 29.2919 & 39.3270& 10.2248 & 59.4 \\
   3 & 29.2824 & 38.5239& 10.4010 & 558.5 \\
   4 & - & 37.6402& - & 1298 \\
   5 &   -      & 36.7088  & - & 1653.3 \\
    6 &   -      & 35.7561  & - & 2033.6 \\
    7 &   -      & 35.2567  & - & 2115.2 \\
    8 &   -      & 34.9022  & - & 1706.1 \\
    \textbf{total time} && & $\bm{31.2841}$ & 9469.1  \\

  \bottomrule
\end{tabular}}
 \caption{Iteration cost and computation time for Example \ref{ex3}.}
 \label{ex3_it_cost}
\end{table}

\begin{figure}[htbp]
\center
\includegraphics[width=2.5in]{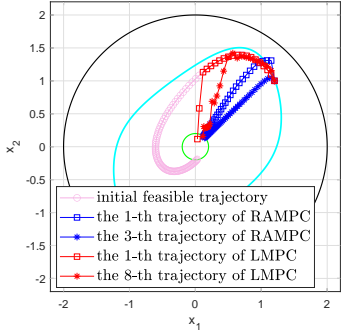}
\caption{Trajectories in Example \ref{ex3} (\textcolor{green}{Green}, black and \textcolor{cyan}{cyan} curve- \textcolor{green}{$\partial \mathcal{T}$},  \textcolor{black}{$\partial \mathcal{X}$} and \textcolor{cyan}{$\partial \mathcal{R}^0$}).}
\label{ex3_figure}
\end{figure}

\end{example}

\begin{example}
\label{ex4}
Consider the Euler version with the time step $\Delta t = 0.1$ of the controlled 3D Van der Pol oscillator from \cite{korda2014controller}:
\begin{equation*}
    \label{3d_vanderpol_sys}
    \begin{cases}
x_1(t+1)=&x_1(t)+\Delta t(-2 x_2(t)) \\
 x_2(t+1)=&x_2(t)+\Delta t(0.8 x_1(t)-2.1 x_2(t)\\
                &+x_3(t)+10 x_1^2(t) x_2(t)) \\
x_3(t+1)=&x_3(t)+\Delta t(-x_3(t)+x^3_3(t))+u(t)
    \end{cases}
\end{equation*}
with the safe set $\mathcal{X}=\{(x_1,x_2,x_3)^{\top}\mid  x_1^2+x_2^2+x_3^2-0.5^2\leq 0\}$, initial state $x_0 = (0.2,0.4,0.1)^\top$, target set $\mathcal{T}=\{(x_1,x_2,x_3)^{\top}\mid x_1^2+x_2^2+x_3^2-0.1^2\leq 0\}$ and input set $\mathcal{U} = \{u \mid -2 \leq u \leq 2\}$. 

The set 
$\mathcal{Y}=\{(x_1,x_2,x_3)^{\top}\mid  x_1^2+x_2^2+x_3^2-0.5\leq 0\}$ is utilized in \eqref{gbf_c} by solving a semi-define program, as described in \cite{zhao2022inner}. We use the cost function $h\left(x,u\right)=\left\|x\right\|_2^2+\left\|u\right\|_2^2$ in this example, where $x = (x_1,x_2,x_3)^\top$. The initial trajectory is generated by the controller $\hat{u}^0(\bm{x})\equiv 0$.  The template of the polynomial form $Q^j_a(\bm{c},x) = c_{1}+c_2x_1+c_3x_2+c_4x_3+c_5x_1x_2+c_6x_1x_3+c_7x_2x_3+c_8x_1^2+c_9x_2^2+c_{10}x_3^2$ and sample number $N'=267$ computed through Theorem \ref{sample_number} are adopted for scenario optimization \eqref{co}. The prediction horizon $N$ and termination condition $\xi$ are same for LMPC. 

In this example, LMPC is unable to achieve the desired task of reducing the iteration cost and generating a suitable controller. In contrast, our RAMPC algorithm demonstrates good performance and can generate an effective controller. Table \ref{ex4_cost_table} summarizes the iteration costs and computation times at each iteration of our RAMPC algorithm. The initial trajectory and ones generated in the first and last iterations of our RAMPC algorithm are visualized in Figure \ref{ex4p_figure}.

\begin{table}[htbp]
    \centering
    \setlength{\tabcolsep}{1.5mm}{
    \begin{tabular}{*{3}{c}}
  \toprule
   \textbf{Iteration}&
   \textbf{Iteration Cost} &\textbf{Time Cost}(seconds)\\

  \midrule
   0 & 1.3489 & \\
   1 & 0.8344 & 7.8186\\
   2 & 0.8295 & 7.8022\\
   3 & 0.8291 & 7.7416\\
   \textbf{total time} & & 23.3623  \\
  \bottomrule
\end{tabular}}
    \caption{Iteration cost and computation time of RAMPC algorithm for Example \ref{ex4}.}
    \label{ex4_cost_table}
\end{table}

\begin{figure}[htbp]
\center
\includegraphics[width=2.5in]{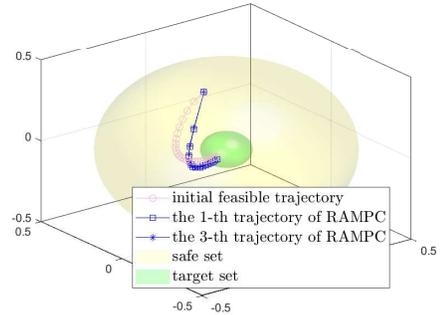} 
\caption{Trajectories in Example \ref{ex4}}
\label{ex4p_figure}
\end{figure}

\end{example}

\paragraph{Discussions on configuration parameters.} Here we give a brief discussion on the configuration parameters: $\lambda,M,K, \xi,N,\delta,\epsilon$ based on the experiments (some are presented in the appendix). The less conservative the computed reach-avoid set in each iteration is, the smaller the iteration cost is. Therefore, $\lambda$ is recommended to be as close to one as possible, but it cannot be one. For details, please refer to \cite{zhao2022inner}; the upper bound $M$ can be any positive number since if $v(x)$ satisfies constraint \eqref{gbf_c} with $M=M'$, then $\frac{v(x)}{M'}$ also satisfies constraint \eqref{gbf_c} with $M=1$; $K$ and $\xi$, especially $\xi$, affect the number of iterations and thus the quality of controllers. Generally, a larger $K$ and/or smaller $\xi$ will improve the control performance, but the computational burden increases; the approximation error for terminal costs is formally characterized in the PAC sense using $\epsilon$ and $\delta$. These two parameters and the size of parameters $\bm{c}$ determine the minimal number of samples for approximating terminal costs according to Theorem \ref{sample_number}. Smaller $\epsilon$ and $\delta$ will lead to more accurate approximations of terminal costs in the PAC sense. However, regarding the overfitting issue and the computational burden of solving \eqref{co}, too small $\epsilon$ and $\delta$ are not desired in practice; the prediction horizon $N$ controls how far into the future the MPC predicts the system response. Like other MPC schemes, a longer horizon generally increases the control performance, but requires an increasingly powerful computing platform (due to solving \eqref{MPC10} online), excluding certain control applications. In practice, an appropriate choice strongly depends on computational resources and real-time requirements. In addition, it is worth remarking here that the costs obtained by RAMPC are not sensitive to $N$ for some cases. For details, please refer to experimental results in the appendix.
\section{Conclusion}
\label{sec:con}
In this paper we proposed a novel RAMPC algorithm for solving a reach-avoid optimization problem. The RAMPC algorithm is built upon MPC and reach-avoid analysis. Rather than solving computationally intractable mixed-integer nonlinear optimization online, it addresses computationally more tractable nonlinear optimization. Moreover, due to the incorporation of reach-avoid analysis, which expands the viability space, the RAMPC algorithm can provide better controllers with a smaller number of iterations. Several numerical examples demonstrated the advantages of our RAMPC algorithm.

\section*{Acknowledgments}
This research was supported by the NSFC under grant No.\ 61836005, CAS Project for Young Scientists in Basic Research under grant No.\ YSBR-040,  NSFC under grant No. 52272329, and CAS Pioneer Hundred Talents Program.
\bibliographystyle{named}
\bibliography{ijcai23}
\clearpage

\section*{Appendix I}
\textbf{The proof of Theorem \ref{feasibility}:}
\begin{proof}

    We firstly prove the RAMPC \eqref{MPC10} satisfies  recursive feasibility for $0\leq t \leq L^j-1$ at every iteration $1\leq j \leq K$.

 %We begin with $j=1$ and then use the induction principle to obtain that the RAMPC \eqref{MPC10} is feasible for $j>1$.

 %Since a feasible state-control trajectory $\{(x^0(i),u^0(i))\}_{i=0}^{L^0}$ is assumed, $\{(x^0(i),u^0(i))\}_{i=0}^{N}$ is a feasible solution to Optimization \eqref{MPC10} with $t=0$. Therefore, \eqref{MPC10} is feasible. Assume that 

 %   \begin{itemize}
 %       \item []\textbf{Step $1$}: at time $t=0$ of the $j$-th iteration the RAMPC \eqref{MPC10} is feasible.
       
%By assumption that the trajectory $\mathbf{x}^0$ is feasible and satisfies the reach-avoid property, At time $t=0$ of the $j$-th iteration the $N$ steps trajectory
 %   \begin{equation}
 %   \label{feasible_x}
 %       \mathbf{x}^0_{0:N} = [x^0_0, x^0_1, \dots, x^0_N] \subset \mathcal{R}^0 
 %   \end{equation}
 %   with the related input sequence,
 %   \begin{equation}
 %   \label{feasible_u}
 %       [\hat{u}^0(x^0_0), \hat{u}^0(x^0_1), \dots, \hat{u}^0(x^0_{N-1})]
 %   \end{equation}
 %   satisfies the constraints in \eqref{MPC10}, thus \eqref{feasible_x} and \eqref{feasible_u} is a feasible solution to RAMPC \eqref{MPC10} at $t=0$ of the $1$-th iteration.

%    Similarly, replace the iteration index $0$ with $1$ in \eqref{feasible_x} and \eqref{feasible_u}, then RAMPC \eqref{MPC10} also has a feasible solution at $t=0$ of the $2$-th{\textcolor{red}{?}} iteration. Repeat this progress then we conclude that at time $t=0$ of the $j$-th iteration the RAMPC \eqref{MPC10} is feasible.

 %    \item []\textbf{Step $2$}: 
1) Assume that the $(j-1)$-th iteration the trajectory generated by RAMPC \eqref{MPC10} is feasible (when $j=1$, the initial trajectory is feasible), then we will show that at time $t=0$ of the $j$-th iteration the RAMPC \eqref{MPC10} is feasible.

At time $t=0$ of the $j$-th iteration the first $N$ steps trajectory
    \begin{equation}
    \label{feasible_x}
    \small
    \{x^{j-1}_0, x^{j-1}_1, \dots, x^{j-1}_N  \} \subset \mathcal{R}^{j-1}
    \end{equation}
    with the related input sequence,
    \begin{equation}
    \label{feasible_u}
    \small
        \{\hat{u}^{j-1}(x^{j-1}_0), \hat{u}^{j-1}(x^{j-1}_1), \dots,\hat{u}^{j-1}(x^{j-1}_{N-1})\}
    \end{equation}
    satisfies the constraints in \eqref{MPC10}, thus \eqref{feasible_x} and \eqref{feasible_u} is a feasible solution of RAMPC \eqref{MPC10} at $t=0$ of the $j$-th iteration.

2).   Assume that at time $t$ of the $j$-th iteration the RAMPC \eqref{MPC10} is feasible, then we will show that at time $t+1$ of the $j$-th iteration RAMPC \eqref{MPC10} is feasible.

    Let $\mathbf{x}_{0:N \mid t}^{*,j}$ and $\mathbf{u}_{0:N \mid t}^{*,j}$  in \eqref{optimal_solution} be the optimal trajectory and input sequence at time $t$. 
%      Note that from constraints in \eqref{MPC10} and \eqref{x_deploy}, the realized state and input at time $t$ of the $j$-th iteration are given by 
%     \begin{equation}
%         \label{realized x_u}
%         \begin{aligned}
% x_t^j & =x_{0 \mid t}^{*, j}, \\
% u_t^j & =u_{0 \mid t}^{*, j} .
% \end{aligned}
%     \end{equation}
    % The terminal constraint \eqref{RA_constraint} enforces $x_{N \mid t}^{*, j} \in \mathcal{R}^1$, 
    
    We define $x_{N\mid t+1}^j \triangleq f(x_{N \mid t}^{*, j},\hat{u}^{j-1}(x_{N \mid t}^{*, j}))$. With the fact that $\mathcal{R}^{j-1}$ is a reach-avoid set with the controller $\hat{u}^{j-1}$, and the first constraint in \eqref{gbf_c}, we have $v(f(x_{N\mid t+1}^j)\geq \lambda v(x_{N \mid t}^{*, j})$ which satisfies the terminate constraint \eqref{terminate_constraint}.
    
    At time $t+1$ of the $j$-th iteration the input sequence 
    \begin{equation}
    \label{feasible_u_t}
    \small
        \{u_{1 \mid t}^{*, j}, u_{2 \mid t}^{*, j}, \dots, u_{N-1 \mid t}^{*, j}, \hat{u}^{j-1}(x_{N \mid t}^{*, j})\}
    \end{equation}
    with the related feasible state trajectory,
    \begin{equation}
    \label{feasible_x_t}
    \small
    \{x_{1 \mid t}^{*, j}, x_{2 \mid t}^{*, j}, \dots, x_{N \mid t}^{*, j},x_{N\mid t+1}^j\}
    \end{equation}
    satisfies the constraints in \eqref{MPC10}, Thus \eqref{feasible_u_t} and \eqref{feasible_x_t} is a feasible solution for RAMPC \eqref{MPC10} at time $t+1$.
     
 %   \end{itemize}

    % The facts has been proven at the $j$-iteration, $\forall j \geq 0$:
    % \begin{itemize}
    %     \item[(1)] The RAMPC is feasible at time $t=0$.
    %     \item[(2)] If the RAMPC is feasible at time $t$, then the RAMPC is feasible at time $t+1$.
    % \end{itemize}
    With the principle of induction, the RAMPC \eqref{MPC10} is feasible for $1\leq j \leq K$ and $0\leq t \leq L^j-1$.

    Next we will prove that the resulting trajectory $\mathbf{x}^j=(x^j(0),\ldots,x^j(L^j))$ (as shown in \eqref{x_deploy}) at iteration $ 1\leq j\leq K$ satisfies the reach-avoid property.
    According to the terminate constraint \eqref{terminate_constraint}, we have that
    \begin{equation*}
    \small
        v(x_0) \leq \frac{v(x_{N\mid 0})}{\lambda^N} \leq \frac{v(x_{N\mid 1})}{\lambda^{N+1}} \leq \cdots \leq \frac{v(x_{N\mid t})}{\lambda^{N+t}} \leq \cdots
    \end{equation*}
i.e. $v(x_{N\mid t}) \geq \lambda^{N+t}v(x_0)$, If the trajectory never reach the target set $\mathcal{T}$, then $\lim_{t \rightarrow +\infty}v(x_{N\mid t}) = +\infty$, which contracts the fact that $v(\cdot)$ is bounded over the compact set $\mathcal{Y}$. Consequently, there exists $t\in \mathbb{N}$ such that $x_{N\mid t} \in \mathcal{T}$ and thus 
$L^j\leq t+N$. Also, from the constraint $v(x)\leq M, \forall x\in \mathcal{T}$, we have that $t+N\leq \log_{\lambda} \frac{M}{v(x_0)}$(see Remark \ref{finite_horizon}) and thus 
$L^j \leq \log_{\lambda} \frac{M}{v(x_0)}$.
%the trajectory will reach the target set $\mathcal{T}$ in finite time $T' \leq \log_{\lambda} \frac{M}{v(x)}$ (see Remark \ref{finite_horizon}).
\end{proof}
%Next we will prove that the $j$-th iteration cost $J^j(x_0)$ is non-increasing as $j$ increases.
%\begin{theorem}
%    \label{J_non_increasing}
%    Consider system \eqref{system} with controllers obtained by solving \eqref{MPC10}. Then, the iteration cost $J^j(x_0)$ is non-increasing with the iteration index $j$.
%\end{theorem}

\textbf{The proof of Theorem \ref{J_non_increasing}: }
\begin{proof}
    Notice that $J_{t \rightarrow t+N}^{\text {RAMPC}, j}(x_t^j) = J_{0 \rightarrow N}^{\text {RAMPC}, j}(x_t^j)$ in \eqref{MPC10}, we replace $J_{t \rightarrow t+N}^{\text {RAMPC}, j}(x_t^j)$ with $J_{0 \rightarrow N}^{\text {RAMPC}, j}(x_t^j)$ in the following proof, we assume $N\leq L^j-1$. %is less than the steps $T'$ of the trajectory at each iteration. 

    Given the optimal state-control pair in \eqref{optimal_solution}, we have 
    \[ \small J_{0 \rightarrow N}^{\text {RAMPC}, j}(x_t^j) = \min_{u_{k\mid t}}\sum_{k=0}^{N-1}h(x_{k\mid t},u_{k\mid t})+Q^{j-1}(x^j_{N\mid t},\pi_{\hat{u}^{j-1}})\]
    \begin{equation}
    \label{J_decrease_t}
    \small
        \begin{aligned}
%            &=h(x_{0 \mid t}^{*, j}, u_{0 \mid t}^{*, j})+\sum_{k=1}^{N-1}h(x_{k \mid t}^{*, j}, u_{k \mid t}^{*, j}) + Q^{j-1}(x_{N \mid t}^{*, j},\pi_{\hat{u}^{j-1}}) \\
            &=h(x_{0 \mid t}^{*, j}, u_{0 \mid t}^{*, j})+\sum_{k=1}^{N-1}h(x_{k \mid t}^{*, j}, u_{k \mid t}^{*, j}) + h(x_{N \mid t}^{*, j}, \hat{u}(x_{N \mid t}^{*, j})) \\
            &~~~~~~~~~~~~~~~~~~~~~~~~~~~~~~~~~~~~~~~~~~~~~~~~+ Q^{j-1}(f(x_{N \mid t}^{*, j}, \hat{u}(x_{N \mid t}^{*, j})) \\
            &\geq h(x_{0 \mid t}^{*, j}, u_{0 \mid t}^{*, j}) + J_{0 \rightarrow N}^{\text {RAMPC}, j}(x_{1 \mid t}^{*, j})
        \end{aligned}
    \end{equation}

    % The last inequation based on the fact that 
    % $$
    From \eqref{x_deploy} and \eqref{J_decrease_t}, we have that for $ x_t^j \in \mathcal{X}$ and $u_t^j \in \mathcal{U}$,
    \begin{equation}
    \begin{aligned}
    \label{J_decrease_t_final}
        J_{0 \rightarrow N}^{\text {RAMPC}, j}(x_{t+1}^{j}) - J_{0 \rightarrow N}^{\text {RAMPC}, j}(x_{t}^{j}) \leq -h(x_t^j,u_t^j) \leq 0.
    \end{aligned}
    \end{equation}
    
    Next, we show that $J^{j-1}(x_0) \geq J_{0 \rightarrow N}^{\text {RAMPC}, j}(x_{0}^{j})$.

    %By the definition of the iteration cost \eqref{iteration cost}, 
    We have 
    \begin{equation}
    \small
    \label{J^j-1>=J^RAMPC}
    \begin{aligned}
        J^{j-1}(x_0) &= \sum_{k=0}^{L^{j-1}-1} h(x_k^{j-1},u_k^{j-1})\\
        &=\sum_{k=0}^{N-1} h(x_k^{j-1},u_k^{j-1}) + \sum_{k=N}^{L^{j-1}-1} h(x_k^{j-1},u_k^{j-1})\\
        &=\sum_{k=N}^{L^{j-1}-1} h(x_k^{j-1},u_k^{j-1}) + Q^{j-1}(x_N^{j-1})\\
        &\geq \min_{u_0,u_1,\dots,u_{N-1}}\Bigg[\sum_{k=0}^{N-1}h(x_k,u_k) +Q^{j-1}(x_N)\Bigg]\\
        &=J_{0 \rightarrow N}^{\text {RAMPC}, j}(x_{0}^{j})
    \end{aligned}    
    \end{equation}

From \eqref{J_decrease_t_final}, we have 
     \begin{equation}
     \small
     \label{J^RAMPC>=J^j}
    \begin{aligned}
        J_{0 \rightarrow N}^{\text {RAMPC}, j}(x_{0}^{j}) &\geq h(x_0^j,u_0^j) + J_{0 \rightarrow N}^{\text {RAMPC}, j}(x_{1}^{j}) \\
        &\geq h(x_0^j,u_0^j) +h(x_1^j,u_1^j)+ J_{0 \rightarrow N}^{\text {RAMPC}, j}(x_{2}^{j})\\
        &\geq \sum_{k=0}^{L^{j}-N-1}h(x_k^j,u_k^j) + J_{0 \rightarrow N}^{\text {RAMPC}, j}(x_{L^{j}-N}^{j})\\
        &\geq \sum_{k=0}^{L^{j}-1}h(x_k^j,u_k^j)\\
        &=J^j(x_0)
    \end{aligned}    
    \end{equation}

    Finally from \eqref{J^j-1>=J^RAMPC} and \eqref{J^RAMPC>=J^j}, we have
    \begin{equation*}
        J^{j-1}(x_0)\geq J_{0 \rightarrow N}^{\text {RAMPC}, j}(x_{0}^{j}) \geq J^j(x_0)
    \end{equation*}
    therefore the iteration cost is non-increasing.
\end{proof}

%\newpage 
\section*{Appendix II}
\begin{example}
\label{ex5}
Consider the system in Example \ref{ex1}, but with the different configurations in Table \ref{ex5_tab1}.
Table \ref{ex5_it_cost} summarizes the costs obtained at each iteration of our RAMPC algorithm and the LMPC one, and Table \ref{ex5_table_time} shows the computation times of each iteration for these two algorithms. Some trajectories generated  by our RAMPC algorithm and the LMPC algorithm are visualized in Figure \ref{ex5_figure}.
\begin{table}[htbp]
    \centering
    \small
    \begin{tabular}{|c|c|c|c|c|c|c|}
    \hline
         $\lambda$& $M$&$N$&$K$&$\xi$&$\delta$&$\epsilon$ \\\hline
          1.001&1&2&10&0.01&0.1&0.1 \\\hline
    \end{tabular}
    \caption{Configuration Parameters in Algorithm \ref{RAMPC_alg} for Example \ref{ex5}.}
    \label{ex5_tab1}
\end{table}

\begin{table}[htbp]
    \centering
    \small
    \setlength{\tabcolsep}{1.5mm}{
    \begin{tabular}{*{3}{c}}
  \toprule
   \multirow{2}*{\textbf{Iteration}}&
   \multicolumn{2}{c}{\textbf{Iteration Cost}}  \\
  \cmidrule(lr){2-3}
  &RAMPC & LMPC \\
  \midrule
   0 &  369.8267 &  369.8267 \\
   1 &  215.1011 & 217.1582\\
   2 & 215.1025 & 216.0086\\
   3 & - & 215.5209\\
   4 & - & 215.2383\\
   5 &   -      & 215.1722  \\
    6 &   -      & 215.1324  \\
    7 &   -      &  215.1106  \\
    8 &   -      & 215.1003  \\
    9 &   -      &  215.1007 \\
  \bottomrule
\end{tabular}}
 \caption{Iteration costs for Example \ref{ex5}.}
 \label{ex5_it_cost}
\end{table}

\begin{table}[htbp]
\small
    \centering
    \setlength{\tabcolsep}{1.5mm}{
    \begin{tabular}{*{3}{c}}
  \toprule
   \multirow{2}*{\textbf{Iteration}}&
   \multicolumn{2}{c}{\textbf{Time Cost}(seconds)}  \\
  \cmidrule(lr){2-3}
  &RAMPC & LMPC \\
   \midrule
   $\bm{1}$ &  3.3884 & 2.6410 \\
   $\bm{2}$ &  2.7991 & 1.4724 \\
   $\bm{3}$ & - &  1.0895 \\
   $\bm{4}$ & - & 1.0568 \\
    $\bm{5}$ & - & 1.0892\\
    $\bm{6}$ & - & 1.0536 \\
    $\bm{7}$ & - & 1.0428 \\
    $\bm{8}$ & - & 1.1520  \\
    $\bm{9}$ & - & 1.0190  \\
    \textbf{total} & $\bm{ 6.1875}$ & 11.6164 \\
   \bottomrule
\end{tabular}}
\caption{Computation times for Example \ref{ex5}.}
\label{ex5_table_time}
\end{table}
\begin{figure}[htbp]
\center
\includegraphics[width=2.5in]{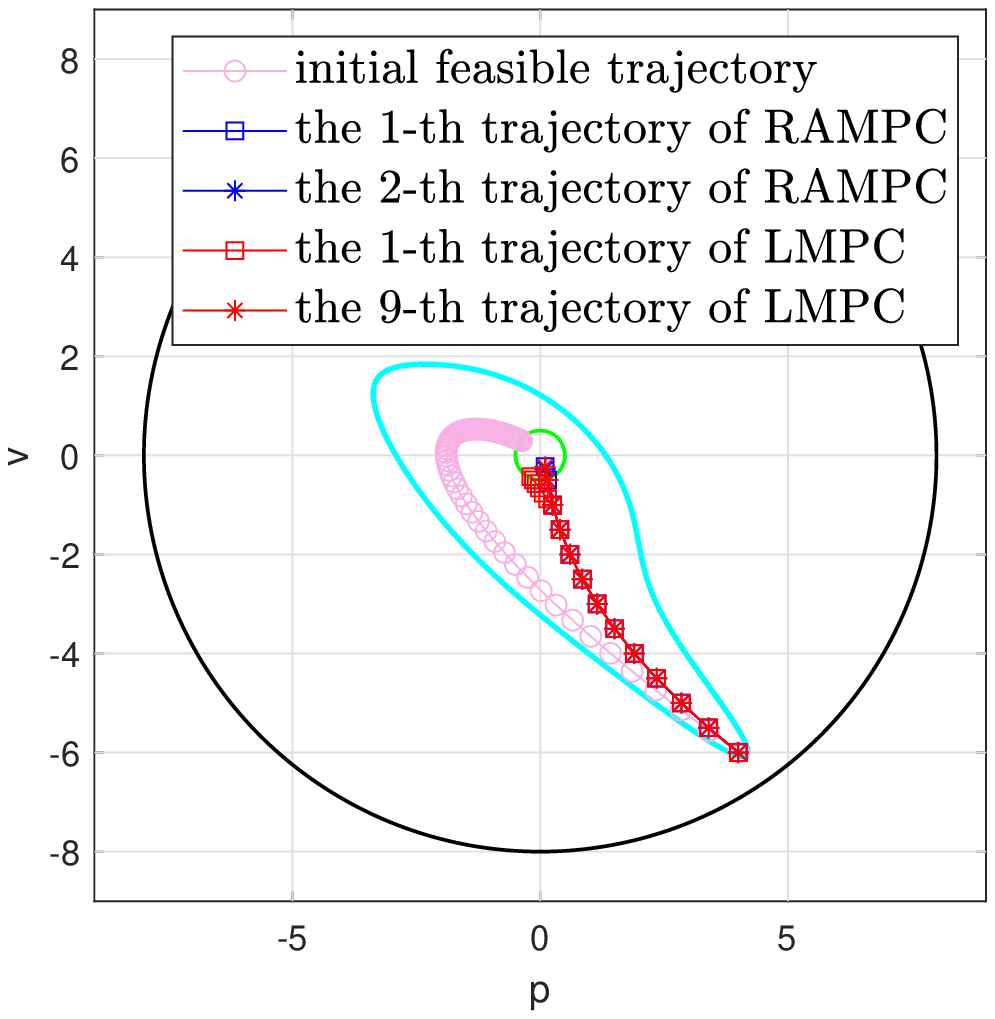}
\caption{Trajectories in Example \ref{ex5}(\textcolor{green}{Green}, black and \textcolor{cyan}{cyan} curve- \textcolor{green}{$\partial \mathcal{T}$},  \textcolor{black}{$\partial \mathcal{X}$} and \textcolor{cyan}{$\partial \mathcal{R}^0$}).}
\label{ex5_figure}
\end{figure}

\end{example}

\begin{table}[htbp]
\small
    \centering
    \begin{tabular}{|c|c|c|c|c|c|c|}
    \hline
         $\lambda$& $M$&$N$&$K$&$\xi$&$\delta$&$\epsilon$ \\\hline
          1.001&1&2&10&0.01&0.05&0.05 \\\hline
    \end{tabular}
    \caption{Configuration Parameters in Algorithm \ref{RAMPC_alg} for Example \ref{ex6}.}
    \label{ex6_tab1}
\end{table}
\begin{example}
\label{ex6}
Consider the system in Example \ref{ex3} with $\Delta t=0.1$ and the different configurations in Table \ref{ex6_tab1}.
Table \ref{ex6_it_cost} summarizes the costs obtained at each iteration of our RAMPC algorithm and the LMPC one, and Table \ref{ex6_table_time} shows the computation times of each iteration for these two algorithms. Some trajectories generated  by our RAMPC algorithm and the LMPC algorithm are visualized in Figure \ref{ex6_figure}. %\textcolor{red}{Since $N=2$ is too small, the LMPC algorithm is unstable and can not reduce the iteration costs after the $1$-th trajectory, probability with the reason the limitation of the terminate state and intractable mixed-integer nonlinear optimization. However, our RAMPC algorithm also works for this case and has a better performance.}

\begin{figure}[htbp]
\center
\includegraphics[width=2.5in]{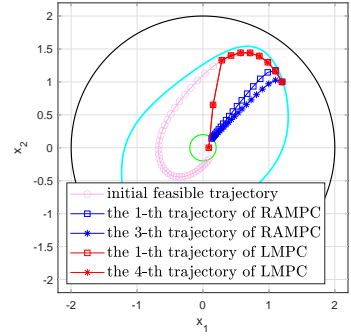}
\caption{Trajectories in Example \ref{ex6} (\textcolor{green}{Green}, black and \textcolor{cyan}{cyan} curve- \textcolor{green}{$\partial \mathcal{T}$},  \textcolor{black}{$\partial \mathcal{X}$} and \textcolor{cyan}{$\partial \mathcal{R}^0$}).}
\label{ex6_figure}
\end{figure}

\begin{table}[htbp]
\small
    \centering
    \setlength{\tabcolsep}{1.5mm}{
    \begin{tabular}{*{3}{c}}
  \toprule
   \multirow{2}*{\textbf{Iteration}}&
   \multicolumn{2}{c}{\textbf{Iteration Cost}}  \\
  \cmidrule(lr){2-3}
  &RAMPC & LMPC \\
  \midrule
   0 &  36.0724 &  36.0724 \\
   1 &  15.4891 & 20.4661\\
   2 & 15.1858 & 21.3654\\
   3 & 15.1860 & 20.4661\\
   4 & - & 20.4661\\
  \bottomrule
\end{tabular}}
 \caption{Iteration costs for Example \ref{ex6}.}
 \label{ex6_it_cost}
\end{table}
\begin{table}[htbp]
    \centering
    \small
    \setlength{\tabcolsep}{1.5mm}{
    \begin{tabular}{*{3}{c}}
  \toprule
   \multirow{2}*{\textbf{Iteration}}&
   \multicolumn{2}{c}{\textbf{Time Cost}(seconds)}  \\
  \cmidrule(lr){2-3}
  &RAMPC & LMPC \\
   \midrule
   $\bm{1}$ &  5.2348  & 208.7463 \\
   $\bm{2}$ & 5.0935 & 142.7198 \\
   $\bm{3}$ &  5.0953 &  9.9840 \\
   $\bm{4}$ & - &  6.4926 \\
    \textbf{total} & $\bm{ 15.4236}$ & 367.9428 \\
   \bottomrule
\end{tabular}}
\caption{Computation times for Example \ref{ex6}.}
\label{ex6_table_time}
\end{table}
\end{example}

\begin{example}
\label{ex2_N}
Consider the system in Example \ref{ex3} with $\Delta t=0.1$: 

\begin{enumerate}
    \item[a] with configuration parameters shown in Table \ref{ex2_N_tab1}.
Table \ref{ex2_N_it_cost} summarizes the costs obtained at each iteration of our RAMPC algorithm with different $N$'s, and Table \ref{ex2_N_table_time} shows the corresponding computation times.
    \item[b] with configuration parameters shown in Table \ref{ex2_N2_tab1}.
Table \ref{ex2_N2_it_cost} summarizes the costs obtained at each iteration of our RAMPC algorithm with different $N$'s, and Table \ref{ex2_N2_table_time} shows the corresponding computation times.
\end{enumerate}

\begin{table*}[htbp]
    \centering
    \begin{tabular}{|c|c|c|c|c|c|c|}
    \hline
         $\lambda$& $M$&$N$&$K$&$\xi$&$\delta$&$\epsilon$ \\\hline
          1.001&1&-&3&0.01&0.05&0.05 \\\hline
    \end{tabular}
    \caption{Configuration Parameters in Algorithm \ref{RAMPC_alg} for Example \ref{ex2_N}a.}
    \label{ex2_N_tab1}
\end{table*}

\begin{table*}[htbp]
    \centering
    \setlength{\tabcolsep}{1.5mm}{
    \begin{tabular}{*{7}{c}}
  \toprule
   \multirow{2}*{\textbf{Iteration}}&
   \multicolumn{6}{c}{\textbf{Iteration Cost}}  \\
  \cmidrule(lr){2-7}
  &2 & 4 & 6 & 8 &10& 12\\
  \midrule
   0 & 36.0724 & 36.0724&36.0724 &  36.0724&  36.0724&36.0724 \\
   1 &  15.4891 &16.4719&15.3693   & 15.3882 &  15.2546&15.1886  \\
   2 &  15.1858 & 15.2198&15.1777  & 15.1640& 15.1755&15.1703  \\
   3 & 15.1860 & 15.1696&15.1673 &  15.1690& 15.1714&15.1719 \\
  \bottomrule
\end{tabular}}
 \caption{Iteration costs with different $N$'s for Example \ref{ex2_N}a.}
 \label{ex2_N_it_cost}
\end{table*}

\begin{table*}[htbp]
    \centering
    \setlength{\tabcolsep}{1.5mm}{
    \begin{tabular}{*{7}{c}}
  \toprule
   \multirow{2}*{\textbf{Iteration}}&
   \multicolumn{6}{c}{\textbf{Time Cost}}  \\
  \cmidrule(lr){2-7}
  &2 & 4 & 6 & 8 &10& 12\\
  \midrule
   1 &  5.2348 &5.2193&  5.9456&3.9578 &4.1903 & 3.5610  \\
   2 &  5.0935  & 5.2343&5.8533 &5.9139 & 6.0883 &6.2036  \\
   3 &  5.0953 & 5.6513&5.6288&6.3013&6.9170&6.4566
  \\
   total & 15.4236 & 16.1050& 17.4276&16.1730& 17.1956&16.2212 \\

  \bottomrule
\end{tabular}}
 \caption{4: Computation times with different $N$'s for Example \ref{ex2_N}a.}
 \label{ex2_N_table_time}
\end{table*}

%\clearpage
%\begin{example}
%\label{ex2_N2}
%Consider the case in Example \ref{ex2}, but 
\begin{table*}[htbp]
    \centering
    \begin{tabular}{|c|c|c|c|c|c|c|}
    \hline
         $\lambda$& $M$&$N$&$K$&$\xi$&$\delta$&$\epsilon$ \\\hline
          1.001&1&-&3&0.01&0.1&0.1 \\\hline
    \end{tabular}
    \caption{Configuration Parameters in Algorithm \ref{RAMPC_alg} for Example \ref{ex2_N}b.}
    \label{ex2_N2_tab1}
\end{table*}

\begin{table*}[htbp]
    \centering
    \setlength{\tabcolsep}{1.5mm}{
    \begin{tabular}{*{7}{c}}
  \toprule
   \multirow{2}*{\textbf{Iteration}}&
   \multicolumn{6}{c}{\textbf{Iteration Cost}}  \\
  \cmidrule(lr){2-7}
  &2 & 4 & 6 & 8 &10& 12\\
  \midrule
   0 & 36.0724 & 36.0724&36.0724 &  36.0724&  36.0724&36.0724 \\
   1 &  17.0569  &16.5073&15.8408   &   15.4489 &  15.2860 & 15.2233   \\
   2 &  15.1790 &15.2286& 15.1805 & 15.1731& 15.1641& 15.1661  \\
   3 & 15.1990 & 15.1678&15.1661 &  15.1777&  15.1922& 15.1703\\
  \bottomrule
\end{tabular}}
 \caption{Iteration costs with different $N$'s for Example \ref{ex2_N}b.}
 \label{ex2_N2_it_cost}
\end{table*}

\begin{table*}[htbp]
    \centering
    \setlength{\tabcolsep}{1.5mm}{
    \begin{tabular}{*{7}{c}}
  \toprule
   \multirow{2}*{\textbf{Iteration}}&
   \multicolumn{6}{c}{\textbf{Time Cost}}  \\
  \cmidrule(lr){2-7}
  &2 & 4 & 6 & 8 &10& 12\\
  \midrule
   1 &  5.4749  &   4.7185 &  5.9308& 6.3891   &  9.3374    & 7.7943    \\
   2 &   5.1461  &  5.5312 &6.5768 & 6.3372   & 6.4644   &8.1602  \\
   3 &   5.4065 &  5.3546&5.6288& 7.1605&  6.8632& 9.3544
  \\
   total & 16.0275
 &  15.6043& 17.4276&  19.8869& 22.6650&  25.3089\\

  \bottomrule
\end{tabular}}
 \caption{Computation times with different $N$'s for Example \ref{ex2_N}b.}
 \label{ex2_N2_table_time}
\end{table*}

\end{example}

\end{document}